\documentclass{amsart}
\usepackage[utf8]{inputenc}
\usepackage{xcolor}
\usepackage{amssymb}
\usepackage{amsmath}
\usepackage{mathtools}
\usepackage{biblatex}

\begin{document}
\begin{titlepage}
   \begin{center}
       \vspace*{1cm}

       \Huge
       \textbf{French to English Mathieu 1873 Translation\\}
       
       \vspace{0.75cm}
       
       \LARGE
       \emph{On the five-fold transitive function of 24 quantities;\\
       Pages 25-46\\}
       
       \vspace{0.75cm}
       
       Original Document Written By \underline{Émile Mathieu}

       \vspace{0.5cm}
       
       \LARGE
        Project Begun on 6th November 2021
            
       \vspace{1.5cm}

       \textbf{Translation Written by}\\
       Yiming Bing, Bright Hu, Ronni Hu, Rhianna Kho, Max Lau, Rhianna Li, Stefan Lu, Finn Mcdonald, Dr Michael Sun, Gavin Trann, Nicholas Wolfe, Joshua Yao, Leon Zhou, Nathan Zhou*
       
       \vspace{1.5cm}
       
       \small
       The link to the original document in French is \underline{http://sites.mathdoc.fr/JMPA/PDF/JMPA\_1873\_2\_18\_A2\_0.pdf}
       \newpage
       
       \Large
       \tableofcontents
            
   \end{center}
\end{titlepage}

In my Memoir \emph{on the functions of several quantities} 
published in volume VI of this Journal, in 1861, I declared (p. 274) that I had a five-fold transitive function of 24 quantities and gave, at the same time, the number of its distinct values.

I had noted two substitutions which leave it invariant, characterizing it completely, so it was easy to verify its existence; but only having this and knowing the number of distinct values sheds no light on the formation of this function. Since then, no one has attempted a proof and I now propose to prove here it actually exists, effectively, and to show how I came to discover it.

I had not published this function in the cited Memoir because the method I was to use lacked the sharpness and elegance of my other results. However, this method will still be of interest to those who can imagine the difficulty of discovering such a function.

\section*{Indication of research process}

{\bf 1.} Let $p$ be a prime number and let $x_0,x_1,x_2,\dots,x_{p-1}$ be $p$ quantities, such that, $x_m = x_n$ whenever $m$ is congruent to $n$ modulo $p$. There are transitive functions of these quantities which are only invariant by substitutions of the form
$$(x_z,x_{az+b});$$
we omit these functions, the symmetric function and the two valued function from further discussions about transitive functions of $p$ letters\footnote{It seems that a transitive function on $p$-letters could somehow correspond to the modern phrase a group of permutations acting transitively on a $p$-element set} because they are well-known.

With this understood, consider a transitive function of a prime number $p$ of quantities, and suppose that $\frac{p-1}{2}$ is also a prime number $q$; I demonstrated in the Memoire cited (Chap. IV), that, by designating the $p$ quantities in a suitable order by $x_0,x_1,\dots,x_{p-1}$ the function is invariant not only by the cyclic substitution

(A) $$(x_0,x_1,x_2,\dots\,x_{p-1}),$$

but also by the substitution

(B) 
$$(x_1,x_{g^2},x_{g^4},\dots,x_{g^{p-3}})(x_{g},x_{g^3},x_{g^5},\dots,x_{g^{p-2}}),$$

where $g$ is a generator modulo $p$; consequently, it is invariant under the $\frac{p(p-1)}{2}$ derived substitutions, included in the formula
$$(z,a^2z+b),$$

$a$ and $b$ being arbitrary.

Now here is the principle on which I relied to determine the function in question.
\\

\emph{Let us represent the substitution (B) by}\\
(C) $(x_0',x_1',x_2',\dots,x_{q-1}')(x_0'',x_1'',x_2'',\dots,x_{q-1}'')$\\\\
\emph{where}
$$x'_0=x_1,x'_1=x_{g^2},x'_2=x_{g^4},\dots,x'_{q-1}=x_{g^p-2};$$ \\
\emph{then, if the second cycle of (B) is similarly identified with the second cycle of (C), the transitive function is necessarily invariant under a regular substitution of $p-3$ quantities of the form} \\
(D) $(x'_z,x'_{y^uz})(x''_z,x''_{y^uz})$

\emph{where $y$ is a generator modulo $q$, $u$ a divisor of $q-1$, smaller than $q-1$, and the indices being taken according to the modulus $q$.}

Before applying this principle to the search for a $23$-letter transitive function, let's apply it to the determination of the $7$- and $11$-letter transitive functions, since $\frac{7-1}{2}$ and $\frac{11-1}{2}$ are prime numbers.

\section*{Transitive functions of 7 and 11 letters}
{\bf 2.} From what we have just recalled, any transitive function of the $7$ quantities $x_0,x_1,x_2,\dots,x_6$ can be considered as invariant under the substitutions
\begin{center}
$(z,z+1)$ and $(z,3^{2}z)$
\end{center}
or by substitutions
$$(e)\quad(x_0,x_1,x_2,x_3,x_4,x_5,x_6),$$
$$(f)\quad(x_1,x_2,x_4)(x_3,x_6,x_5).$$
Let's put these last substitutions in the form
$$(g)\quad(x'_0,x'_1,x'_2)(x''_0,x''_1,x''_2)$$
and we will have the substitution (D)
$$(h)\quad(x'_1,x'_2)(x''_1,x''_2)$$
To identify (f) and (g), we will start by doing\\
\begin{center}
$x'_0=x_1$, $x'_1=x_2$, $x'_2=x_4$,\\
\end{center}
\vspace{0.75cm}
and we will have, by considering the second cycles, the following three hypotheses:\\
\begin{center}
$1^\circ$\quad$x''_0=x_3, x''_1=x_6, x''_2=x_5;$\\
$2^\circ$\quad$x''_0=x_6, x''_1=x_5, x''_2=x_3;$\\ 
$3^\circ$\quad$x''_0=x_5, x''_1=x_3, x''_2=x_6.$\\
\end{center}
\vspace{0.75cm}
And, according to that, the substitution (h) will be one of the following three:\\
\begin{center}
$1^\circ$    $(x_2 x_4) (x_6 x_5),$\\
$2^\circ$    $(x_2 x_4) (x_5 x_3),$\\
$3^\circ$    $(x_2 x_4) (x_3 x_6).$\\
\end{center}
\vspace{0.75cm}
Any of the substitutions $1^\circ,2^\circ,3^\circ$, combined with (e) and (f), gives a system of $4\times6\times7$ conjugated substitutions, which belongs to a twice transitive function which has $\frac{2\cdot3\cdot4\cdot5}{4}=30$ values. These three systems of conjugated substitutions or the three corresponding functions differ only in the way of designating the $7$ quantities.

Consider, for example, the first of these three functions. From the substitutions (f) and $1^\circ$, we immediately see that it is considered as a function of the $6$ quantities $x_1,x_2,x_3,x_4,x_5,x_6$ it is not only invariable by $1^\circ$, but also by
$$(x_2x_6)(x_4x_5);$$
the function is therefore transitive with respect to $x_2,x_4,x_6,x_3$, and it has indeed $30$ values.

According to a general theorem, a transitive function being necessarily invariable by the substitution (f), it suffices, to characterise the current function, to say that it is invariant by the substitutions ($e$) and $1^\circ$.

Let $x_{\infty}$ denote an eighth quantity; there is a function of the $8$ letters $x_0,x_1,\dots,x_6,x_{\infty}$ which, considered as a function of the first $7$ letters, is identical to the previous function, so that it is invariant by the substitutions ($e$) and $1^\circ$; it is, moreover, invariant by all the substitutions
$$\left(z,\frac{Az+B}{Cz+D}\right),$$ \\
for which $AD$-$BC$ is a quadratic residue of $7$; it will be characterised by the substitutions ($e$) and $1^\circ$, and any of the latter substitutions containing $x_{\infty}$, for example the following:

$$(x_zx_{-\frac{1}{z}})=(x_0x_{\infty})(x_1x_6)(x_2x_3)(x_4x_5).$$

This three times transitive function of $8$ quantities is enclosed in a family of twice transitive functions of $p^v$ quantities, which become three times transitive for $p=2$ and which have $\frac{1\cdot2\dots{p^v}}{p^v(p^v-1)\dots(p^v-p^{v-1})}$ values. I have given this family of functions in Chapter III of the Memoir cited, and I have given it a very elegant method of forming them. I have also shown that these functions should be used as solvers of equations whose degree is a power of a prime number. (See \emph{Memoir on the resolution of equations}, in the \emph{Annals of Pure and Applied Mathematics}, t. IV, 1862.)

{\bf 3.} Let us apply the same method to the determination of the transitive functions of $11$ letters.

First, any transitive function of $11$ quantities can be considered as invariant under the substitutions

$$(z,z+1),\qquad(z,2^{2}z),$$
or by substitutions
$$(k)\quad(x_0,x_1,x_2,x_3,x_4,x_5,x_6,x_7,x_8,x_9,x_{10}),$$
$$(l)\quad(x_1,x_4,x_5,x_9,x_3)(x_2,x_8,x_{10},x_7,x_6).$$
\vspace{0.75cm}
Let us represent the last substitution by
\begin{center}
$(m)\quad(x'_0,x'_1,x'_2,x'_3,x'_4)(x''_0,x''_1,x''_2,x''_3,x''_4)$
\end{center}
\vspace{0.75cm}
identify ($m$) with ($l$) starting by putting

$$x'_0=x_1,x'_1=x_4,x'_2=x_5,x'_3=x_9,x'_4=x_3,$$

and we will have to make one of the following five assumptions:\\
\begin{center}
$1^\circ$ $x''_0=x_2, \: x''_1=x_8, x''_2=x_{10}, x''_3=x_7, x''_4=x_6;$\\
$2^\circ$ $x''_0=x_8, \: x''_1=x_{10}, x''_2=x_7, x''_3=x_6, \: x''_4=x_2;$\\
$3^\circ$ $x''_0=x_{10}, x''_1=x_7, x''_2=x_6, x''_3=x_2, x''_4=x_8;$\\
$4^\circ$ $x''_0=x_7, \: x''_1=x_6, x''_2=x_2, x''_3=x_8, x''_4=x_{10};$\\
$5^\circ$ $x''_0=x_6, \: x''_1=x_2, x''_2=x_8, x''_3=x_{10}, x''_4=x_7.$
\end{center}
\vspace{0.75cm}
The substitution
$$(D)\quad(x'_zx'_{y^uz})(x''_zx''_{y^uz}),$$
whose indices are taken according to modulus $5$, will be, taking $y=2$, depending on whether $u=2$ or $u=1$,
$$(n)\quad(x'_1,x'_4)(x'_2x'_3)(x''_1x''_4)(x''_2x''_3),$$
or
$$(p)\quad(x'_1x'_2x'_4x'_3)(x''_1x''_2x''_4x''_3)$$.

Substitution ($n$) is the second power of substitution ($p$). So the system of conjugate substitutions deduced from ($k$) and ($p$) obviously contains the system of conjugate substitutions deduced from ($k$) and ($n$); therefore, if the system deduced from ($k$) and ($p$) gives a transitive function of $11$ letters, we are sure that the other system also gives one, while the converse is not obvious.

So let us start by dealing with the conjugate substitution system deduced from $(k)$ and $(n)$. Based on the five possible assumptions, $(n)$ is one of the following five substitutions:

\begin{center}
$1^\circ$    $(x_4 x_3)(x_5 x_9)(x_8 x_6)(x_{10} x_7),$\\
$2^\circ$    $(x_4 x_3)(x_5 x_9)(x_{10} x_2)(x_7 x_6),$\\
$3^\circ$    $(x_4 x_3)(x_5 x_9)(x_7 x_8)(x_6 x_2),$\\
$4^\circ$    $(x_4 x_3)(x_5 x_9)(x_6 x_{10})(x_2 x_8),$\\
$5^\circ$    $(x_4 x_3)(x_5 x_9)(x_2 x_7)(x_8 x_{10}).$\\
\end{center}

The substitutions $1^\circ,3^\circ,5^\circ$ are not suitable, that is to say that the system of conjugate substitutions derived from $(k)$ and from one of these three substitutions is nothing other than the $\frac{1\cdot2\cdot3\dots11}{2}$ substitutions which leave the function of the $11$ letters which has two values only. We have to show how we can recognize it.

This is very easy if we are willing to make use of this principle, which has served me many times in my research, but not in my demonstrations: it is that a transitive function of a prime number $p$ quantities cannot be invariant by a substitution taking place on less than $\frac{p+1}{2}$ quantities.\footnote{In a Memoir presented to the Academy of Sciences on May 31, 1858, I demonstrated the following theorem, which could also serve the goal proposed here: \emph{If a transitive function of n letters is transitive with respect to a certain number a of its letters which is not greater than $\frac{n}{2}$, these n letters can be divided into groups each composed of the same number of letters, so that the function is transitive with respect to these groups, which are transitive.} But, although I have never doubted the correctness of this theorem, the proof not having fully satisfied me, I have never forgotten it.}

Consider, for example, the system deduced from $(k)$ and $1^\circ$. Let us do $1^\circ$ once, then $(k)$ twice, and form the table.

\begin{center}
$0,1,2,3,4,5,6,7,8,9,10$\\
$0,1,2,4,3,9,8,10,6,5,7$\\
$2,3,4,6,5,0,10,1,8,7,9.$\\
\end{center}

We will have made, ultimately, the substitution $(0,2,4,5)(1,3,5,10,9,7)$, which is part of the system of conjugate substitutions; do this last substitution six times, we will have the substitution
$$(q)\quad(0,4)(2,5),$$
which is also part of this system, and which has only four letters; therefore, according to the above principle, substitution $1^\circ$ should be rejected.

If we do not want to use this principle, here is how we can note that there are no transitive functions invariable by (k) and $1^\circ$. We will deduce from these two substitutions some others which do not not contain the same letter, $x_0$ for example. Most often, we will be able to be satisfied with two substitutions which do not contain $x_0$; we will notice that the function is transitive with respect to the letters remaining $x_1,x_2,\dots,x_{10}$. From the substitutions obtained on these $10$ quantities, let us deduce some others not containing $x_1$, and note that the function, according to these last substitutions, is transitive with respect to the $9$ letters $x_2,x_3,\dots,x_{10}$. We then determine derived substitutions not containing $x_2$, and so on. Continuing thus, we would prove that the function is transitive with respect to $10$ letters, then to $9$ takes from the previous ones, then to $8$ takes from the latter, then to $7$, finally to $6$ letters.

However, a function of $n$ letters cannot be more than $\frac{n}{2}$ times transitive when it has more than two values (see our Memoir \emph{On the Number of Values that a Function can Acquire}, t. V of this Journal, 1860, p. 18); therefore the considered function cannot have more than two values.\\



\indent This method is very long, but it can be simplified a lot. Indeed, in the current case for example, the function is invariant by the substitution (q). So, as soon as we recognize that the function is four times transitive, we can conclude that it is transitive with respect to any $4$ letters, and that it has only two values, since instead of the $4$ letters of the substitution $(q)$ we can bring any $4$ letters into the function.\\
\indent It will thus be recognized that substitution $1^\circ$ must be rejected, and it is the same for substitutions $3^\circ$ and $5^\circ$. It remains to consider the two substitutions $2^\circ$ and $4^\circ$ which are suitable.\\
\indent The function characterized by $(k)$ and $2^\circ$ is twice transitive, and, considered as a function of the 9 quantities $x_2,x_3,\dots,x_{10}$, it is invariant only by substitutions derived from
$$(r)\quad(2,4,7)(3,10,6)(5,8,9)$$
and substitution $2^\circ$ or
$$(4,3)(5,9)(10,2)(7,6);$$
\vspace{0.75cm}
the substitutions derived from these last two are the second power of $(r)$ and these other two
\begin{center}
$(7,10)(8,5)(6,4)(2,3),$\\
$(2,6)(9,8)(3,7)(4,10).$
\end{center}
\vspace{0.75cm}
So, by swapping the $9$ quantities $x_2,x_3,\dots,x_{10}$, we get $6$ equal values of the function, and this function, which is twice transitive, has $\frac{1\cdot2\cdot3\cdot4\cdot5\cdot6\cdot7\cdot8\cdot9}{6}=60,480$ values. It was first given by M. Kronecker.

The transitive function characterised by (k) and $4^\circ$ is twice transitive, and, considered as a function of the $9$ letters $x_2,x_3,\dots,x_{10}$, it is invariable by the substitution
$$(7,4,3)(2,6,9)(8,5,10),$$ \\
\vspace{0.75cm}
the substitution $4^\circ$ and those derived from it. It is easy to see that this function differs from the preceding one only by the way of designating the $11$ quantities which it contains.

Thus, from the examination of the substitution (n), it follows that \emph{there exists a twice transitive function of eleven letters which has $60480$ values and characterized by the substitutions}
$$(E)\quad(x_0,x_1,x_2,\dots,x_{10}),$$
$$(H)(x_4,x_3)(x_5,x_9)(x_{10},x_2)(x_7,x_6).$$\\
\vspace{0.75cm}
\emph{By adding any of the substitutions}
$$(z, \frac{Az+B}{Cz+D} )$$

\emph{for which AD-BC is a quadratic residue of 11 and containing $x_\infty$, for example by}
$$(z,-\frac{1}{z}) = (\,0,\infty\,)\;(\,1,10)\;(\,2,5)\;(\,3,7)\,(\,4,8)\,(\,6,9\,)$$
\emph{which we have the three substitutions which characterize three times transitive function of the twelve quantities, $x_\infty, x_0, x_1, \dots, x_{10}$}, which, considered as a function of the last eleven quantities, is identical to the twice transitive function which precedes.

{\bf 4}. Now let's consider the substitution (\emph{p}). The first, third and fifth identifications, having to be rejected in the examination of the substitution (\emph{n}), all the more so now; it is therefore necessary to try only the second and the fourth identification, which gives, for the substitution (\emph{p}), one of these two formulas

$$(s)\quad(x_4,x_5,x_3,x_9)(x_{10},x_7,x_2,x_6),$$
$$(t)\quad(x_4,x_5,x_3,x_9)(x_6,x_{10},x_2,x_8),$$
\\
There exists a function transitive of 11 letters, characterised by (\emph{k}) and (\emph{s}) or by (\emph{k}) and (\emph{t}). Since the substitutions (\emph{k}) and (\emph{s}) give the same system of conjugated substitutions as (\emph{k}) and (\emph{t}), let us consider only the first case.\\

The function characterised by (\emph{k}) and (\emph{s}) is evidently invariable by the whole system of conjugate substitutions of the function of Mr Kronecker; as a consequence, it is invairable by the substitution (\emph{r}), and, adding it to (\emph{l}) and (\emph{s}) we have all three substitutions.\\
(L) $(x_1,x_4,x_5,x_9,x_3)(x_2,x_8,x_7,x_6)$\\
(R) $(x_2,x_4,x_7)(x_3,x_10,x_6)(x_5,x_8,x_9)$\\
(S) $(x_4,x_5,x_3,x_9)(x_10,x_7,x_2,x_6)$\\

We can recognize that the system of conjugate substitutions derived from these three substitutions, which are carried out on the $10$ quantities $x_2,x_2,\dots,x_{10}$ is three times transitive of $11$ letters which has $7!$ values

Thus, from the examination of the substitution (p), it follows that \emph{there is a four times transitive function $\gamma$ of eleven letters which has $7!=5040$ values characterised by the substitutions}

$$(x_0,x_1,x_2,\dots,x_{10},)$$
$$(x_4,x_5,x_3,x_9)(x_{10},x_7,x_2,x_6),$$

\emph{in whose system of conjugate substitutions contains the system of substitutions of the twice transitive function characterised by the substitutions (E) and (H). if to these two substitutions we add the substitution $(z, -\frac{1}{z}),$ we have the three substitutions which characterise a five times transitive function of twelve quantities}.

I present the discovery of this five-fold transitive function in a very different way at the end of Chapter II of my Memoir \emph{on the functions of several quantities}. It is useful to compare the form in which it is presented here to the one in which we first made it known.

Let $\omega$ be a root of an irreducible quadratic taken with respect to modulus 3, for example a root of the congruence

$$\omega^2 + 2\omega + 2 \equiv 0 (\mathrm{mod}\quad3);$$
all the roots of

$$x^{3^2} \equiv x(\mathrm{mod}\quad3)$$
 are equal to

$$0,1,2,\omega,1+\omega,2+\omega,2\omega,1+2\omega,2+2\omega;$$
the quantity chosen for $\omega$ is, moreover, a primitive root of the congruence

$$x^{3^2-1}\equiv1(\mathrm{mod}\quad3);$$


\noindent and we will distinguish the 8 roots of this congruence into two categories: the first, which contains the odd powers of $\omega$,

$$\omega, \omega^2\equiv1+2\omega,\quad\omega^5\equiv1+2\omega,\quad\omega^7\equiv1+\omega $$
and which we call \emph{quadratic non-residues}; the second, which contains the even powers of $\omega$,

$$\omega^2\equiv1+\omega,\omega^4\equiv2,\omega^6\equiv2+2\omega,\omega^8\equiv1$$
which we call \emph{quadratic residues}.

According to this, let us change notations and denote the ten quantities

$$x_2, x_3, x_{4}, x_{5}, x_{6}, x_{7}, x_{9}, x_{10}, x_{8}, x_{1}$$
repectively by

$$\gamma_\omega, \gamma_{1+\omega}, \gamma_{2+2\omega}, \gamma_{1+2\gamma}, \gamma_2, \gamma_1, \gamma_{2+\omega}, \gamma_{2\omega}, \gamma_0, \gamma_\infty$$
then the three substitutions (L), (R), (S) become:

$$(\gamma_\infty, \gamma_{2+2\omega}, \gamma_{1+2\omega}, \gamma_{2+\omega}, \gamma_{1+\omega})(\gamma_\omega, \gamma_0, \gamma_{2\omega}, \gamma_1, \gamma_2),$$
$$\gamma_\omega,\gamma_{2+2\omega}, \gamma_1)(\gamma_{1+\omega}, \gamma_{2\omega}, \gamma_2)(\gamma_{1+2\omega}, \gamma_0, \gamma_{2+\omega}),$$
$$(\gamma_{2+2\omega}, \gamma_{1+2\omega}, \gamma_{1+\omega}, \gamma_{2+\omega})(\gamma_{2\omega}, \gamma_1, \gamma_\omega, \gamma_2),$$

\noindent and can be put in this form
$$(z,\frac{z+2\omega}{(1+\omega)z+1}),$$
$$(z,z+2+\omega),$$
$$(z,\omega{z^3});$$

We can easily conclude that all the substitutions of 8, 9 and 10 of the quantities represented by the notion $\gamma$ are included in the two formulas\\

$$(\gamma_z,\gamma_{\frac{Az+B}{Cz+D}}),(\gamma_z,\gamma_{\frac{A_1z^3+B_1}{C_1z^3+D_1}}),$$


\noindent where $AD-BC$ is quadratic residue and $A_1D_1-B_1C_1$ non-residue, and these two formulas represent a system of $8\dot9\dot10$ conjugated substitutions, three times transitive.

We have relied on this theorem that any transitive function of eleven quantities is invariant by the substitution (l) or
$$(z, 2^2z);$$

\noindent but this theorem does not prove that such a function cannot be invariant by the substitution $(z,2z)$ or\\
(u) $$(x_1x_2x_4x_8x_5x_{10}x_9x_7x_3x_6),$$ 

of which the preceding is a second power. To recognize if such a function can exist, we will combine this substitution with the two substitutions (E) and (H), which characterise the twice transitive function of 11 letters, and we will see that the system of conjugate substitutes which s' would deduce from it would be none other than the $1.2.3\dots11$ substitutions which leave the symmetric function invariant. It follows that there are no transitive functions of $11$ letters invariant by the substitution (u), or more generally by no even substitution. (We call \emph{pair substitution} the one that changes the function that has two values.)

\section*{Functions acting transitively on $23$ elements}

\textbf{5.} Since the number $2$3 is a prime number, we can apply the same method to the determination of the transitive functions of $23$ quantities. 

From what we know, any transitive function of $23$ quantities $x_0, x_1, x_2, \cdots, x_{22}$ can be considered invariant by the substitutions
$$(z, z+1) \; \text{and} \; (z, 5^2z),$$

since $5$ is the primitive root of $23$, or by substitutions\\
\\
(A)$(x_0,x_1,x_2,x_3,\dots,x_{22})$\\
$$\text{(B)}\begin{cases}
(x_1,x_2,x_4,x_8,x_{16},x_9, x_{18},x_{13},x_3,x_6,x_{12}),\\
(x_5,x_{10},x_{20},x_{17},x_{11},x_{22},x_{21},x_{19},x_{15},x_7,x_{14}).\\
\end{cases}$$\\
\\
Let's put this last substitution in this form

$$\text{(C)}\begin{cases}
(x^\prime_0,x'_1,x'_2,x'_3,x'_4,x'_5,x_6',x'_7,x'_8,x'_9,x'_{10})\\
(x^{\prime\prime}_0,x''_1,x''_2,x''_3,x''_4,x''_5,x''_6,x''_7,x''_8,x''_9,x''_{10})
\end{cases}$$\\

and let $x'_0,x'_1,\dots,x'_{10}$ be equal to the quantities of the same row of the first cycle of (B), so that we have
$$x'_0=x_1,x'_1=x_2,x'_2=x_4,x'_3=x_8,x'_4=x_{16},x'_5=x_9,x'_6=x_{18},x'_7=x_{13},x'_8=x_3,x'_9=x_6,x'_{10}=x_{12};$$

but the identification of the second cycle of (C) with the second cycle of (B) can be done in eleven different ways:\\
$$1^\circ\begin{cases}
x''_0=x_5,x''_1=x_{10},x''_2=x_{20},x''_3=x_{17},x''_4=x_{11},x''_5=x_{22},\\
x''_6=x_{21},x''_7=x_{19},x''_8=x_{15},x''_9=x_7,x''_{10}=x_{14};
\end{cases}$$
$$2^\circ\begin{cases}
x''_0=x_{10},x''_1=x_{20},x''_2=x_{17},x''_3=x_{11},x''_4=x_{22},x''_5=x_{21},\\
x''_6=x_{19},x''_7=x_{15},x''_8=x_7,x''_9=x_{14},x''_{10}=x_5 

\end{cases}$$

...............................................................................................................................

For the substitution (D) of n°1, we have to take the substitution

$$(T)\begin{cases}
(x'_1,x'_2,x'_4,x'_8,x'_5,x'_{10},x'_9,x'_7,x'_3,x'_6,)\\
(x''_1,x''_2,x''_4,x''_8,x''_5,x''_{10},x''_9,x''_7,x''_3,x''_6,)
\end{cases}$$
composed of two cycles of 10 quantities, or its fifth power, or its second power. Its fifth power is
$$(x'_1x'_{10})(x'_2x'_9)(x'_3x'_8)(x'_4x'_7)(x'_5x'_6),$$
$$(x''_1x''_{10})(x''_2x''_9)(x''_3x''_8)(x''_4x''_7)(x''_5x''_6),$$

and, after having applied one of the $11$ identifications, we can verify that it must be rejected. A fortiori, the substitution (T) must be set aside.

It therefore remains to consider the second power of (T), which is
$$(x_1'x_4'x_5'x_9'x_3')(x_2'x_8'x_{10}'x_7'x_6'),$$
$$(x_1''x_4''x_5''x_9''x_3'')(x_2''x_8''x_{10}''x_7''x_6''),$$
If we apply our original identifications we get the substitution
$$\text{(U)}\begin{cases}
(x_2,x_{16},x_9,x_5,x_8)(x_4,x_3,x_{12},x_{13},x_{18})\\
(x_{10},x_{11},x_{22},x_7,x_{17})(x_{20},x_{15},x_{14},x_{19},x_{21})
\end{cases}$$

which provides a function that we will examine. As for the 10 other identifications, they won't give new functions. We take functions generated by $(A) and (U)$ which omits $x_0$ and hence act transitively on $22$ quantities, namely $x_1, x_2,x_3,\dots,x_{22}$; we will then form substitutions which, besides $x_0$, also omits $x_1$ and get a transitive action on 21 quantities $x_2,x_3,x_4,\dots, x_{22}$. We can also omit $x_0,x_1,x_5$ (from (U)) to get functions acting transitively on $20$ quantities.
These functions therefore act four fold transitively. 

Finally, consider the substitutions which do not contain $x_0$, $x_1$, $x_5$, and a fourth quantity, for example $x_2$. We derive: $1^\circ$ substitution
$$\text{(Z)}\begin{cases}
(4,19)(15,12)(16,7)(10,8)(18,9)(6,13)(11,20)(21,22),\\
(4,16)(19,7)(15,10)(12,8)(18,11)(9,20)(6,21)(13,22),\\
(4,18)(19,9)(15,6)(12,13)(16,11)(7,20)(10,21)(8,22),\\
(4,12)(19,15)(16,8)(7,10)(18,13)(9,6)(11,22)(20,21),
\end{cases}$$

according to which, as we can see at first glance, the function is transitive with respect to the 16 quantities which it has as indices: $4,6,7,8,9,10,11,12,13,15,16,18,19,20,21.22;$ \\

$2^\circ$ substitution \\
$(3,14,17)(7,19,21)(13,15,12)(4,10,11)(9,8,16)(20,18,22).$ \\

We can conclude that there are $16\times3$ substitutions which do not contain $x_0, x_1, x_5, x_2$ and leave the function invariant, given the four times transitive function of $23$ quantities has $\frac{1 \cdot 2 \cdot 3 \cdot 4 \cdots 19}{16\times3}$ values.
We are now going to make a very curious remark: it is that, if we represent the substitution ($U$) cycle by cycle,by
$$\left(\begin{array}{l}\left(\gamma_{0} \gamma_{1}, \gamma_{2} \gamma_{3} \gamma_{4}\right)\left(\gamma_{0}^{\prime} \gamma_{1}^{\prime} \gamma_{2}^{\prime} \gamma_{3}^{\prime} \gamma_{4}^{\prime}\right) \\ \left(\gamma_{0}^{\prime \prime} \gamma_{1}^{\prime \prime} \gamma_{2}^{\prime \prime} \gamma_{3}^{\prime \prime} \gamma_{4}^{\prime \prime}\right)\left(\gamma_{0}^{\prime \prime \prime} \gamma_{1}^{\prime \prime \prime} \gamma_{2}^{\prime \prime \prime} \gamma_{3}^{\prime \prime \prime} \gamma_{4}^{\prime \prime \prime}\right)\end{array}\right)$$
the function is invariant by the substitution $\left(\gamma_{z} \gamma_{4 z}\right)$ with the indices being taken according to $\mod{5}$, i.e. by
$$\left\{\begin{array}{l}\left(\gamma_{1} \gamma_{4}\right)\left(\gamma_{2} \gamma_{3}\right)\left(\gamma_{1}^{\prime} \gamma_{4}^{\prime}\right)\left(\gamma_{2}^{\prime} \gamma_{3}^{\prime}\right) \\ \left(\gamma_{1}^{\prime \prime} \gamma_{4}^{\prime \prime}\right)\left(\gamma_{2}^{\prime \prime} \gamma_{3}^{\prime \prime}\right)\left(\gamma_{1}^{\prime \prime \prime} \gamma_{4}^{\prime \prime \prime}\right)\left(\gamma_{2}^{\prime \prime \prime} \gamma_{3}^{\prime \prime \prime}\right)\end{array}\right)$$
provided that we correctly identify ($U$) and ($V$). The 20 letters of the substitution ($U$) consist of the 16 letters of the substitutions ($Z$) and, in addition, of $x_2,x_3,x_{17},x_{14}$, which are in four different cycles of ($U$). Then, so that ($P$) can coincide with one of the substitutions ($Z$) or one of their derivatives, which are all similar to them, we will do
$$y_0=x_2,y_0'=x_3,y_0''=x_{17},y_0'''=x_{14}$$
by identifying ($U$) and ($V$), and it will result
$$y_1=x_{16},y_2=x_9,y_3=x_6,y_4=x_8$$
$$y_1'=x_{12},y_2'=x_{13},y_3'=x_{18},y_4'=x_4$$
$$y_1''=x_{10},x_2''=x_{11},y_3''=x_{22},y_4''=x_7$$
$$y_1'''=x_{10},y_2'''=x_{21},y_3'''=x_{20},y_{10}'''=x_{15}$$
We conclude, for the substitution ($P$)\\
$$(H)\begin{cases}
(x_{16}x_8)(x_9x_6)(x_{12}x_4)(x_{12}x_{18})\\
(x_{10}x_7)(x_{11}x_{22})(x_{19}x_{15})(x_{21}x_{20})\\
\end{cases}
$$
which precisely represents the fourth of the substitutions ($Z$).

If we take the notation of $x^{\prime}$ and $x^{\prime}$, this substitution could be written

$$\text{(K)}\left\{\begin{array}{l}
\left(x_{4}^{\prime} x_{3}^{\prime}\right)\left(x_{5}^{\prime} x_{0}^{\prime}\right)\left(x_{10}^{\prime} x_{2}^{\prime}\right)\left(x_{7}^{\prime} x_{6}^{\prime}\right) \\
\left(x_{1}^{\prime}, x_{9}^{v}\right)\left(x_{4}^{\prime} x_{5}^{\prime}\right)\left(x_{7}^{\prime} x_{8}^{*}\right)\left(x_{8}^{\prime} x_{2}^{*}\right)
  \end{array}\right.$$

Now the first cycle of $(B)$ or of $(C)$ and the first line of $(K)$ are the two substitutions
($L$)
$$
\begin{array}{l}
\left(x_{0}^{\prime} x_{1}^{\prime} x_{2}^{\prime} \ldots x_{10}^{\prime}\right) \\
\left(x_{9}^{\prime} x_{3}^{\prime}\right)\left(x_{\mathrm{b}}^{\prime} x_{\theta}^{\prime}\right)\left(x_{10}^{\prime} x_{t}^{\prime}\right)\left(x_{z}^{\prime} x_{\theta}^{\prime}\right)
\end{array}
$$

which characterise a twice transitive function of 11 quantities which has $6\times10\times{11}$ equal values (see $\mathbf{n}^{\circ}\mathbf{5}$). The second cycle of ($C$) and the second line of ($K$) give the two substitutions

$$\text{(L')} x_6'x_1^nx_2''\dots{x_{10}'}$$
$$\text{(M')} (x_1''x_9'')(x_4''x_5'')(x_7''x_8'')(x_6''x_2'')$$

which also characterise a twice transitive function of 11 quantities.

\emph{So the four times transitive function of $23$ quantities remains invariable if we make on the $11$ quantities of the first cycle of ($B$) the $6\times10\times1$ I substitutions derived from ($1$) and ($M$) which characterise the twice transitive function of these $11$ quantities, provided that at the same time we do on the quantities of the second cycle of ($B$) sub. Similar statements similarly derived from ($L'$) and ($M'$), which also characterize a twice transitive function.}

\emph{There exists a five times transitive function of $24$ quantities, which, considered as a function of the first $23$ $x_{0},x_{1},\dots,x_{22}$, is identical to the function we have just obtained.}

Let us denote the twenty-fourth letter by $x_{\infty}$; this function is invariant by the substitutions
$$\left(z, \frac{A z+B}{C z+D}\right)$$

for which $A D-B C$ is a quadratic residue; therefore, it is characterized by two substitutions which characterize the function 
precedent of $23$ quantities and not any of the last substitutions containing $x_{\infty}$, for example by $(z,-\frac{1}{z})$,

$$ (0,\infty) (1,22) (2,11) (3,15) (4,17) (5,9) $$
$$  (6,19) (7,13) (8,20) (10,16) (12,21) (18,14)  $$

From what we have seen, any transitive function of $23$ quantities can be considered invariant by substitution ($A$) and the substitution ($U$); but it is still necessary to wonder if, besides the function which we obtain, there does not exist another one that is invariant not only by ($z$, $5^2$z), but also by the circular substitution of 22 letters ($z$, $5z$). However, it will easily be recognized that such a function does not exist.

Now that there is only one function transitive of $23$ quantities, it is evident that that function is characterised by substitution and by any substitutions which leave it invariant, provided that the latter can be taken outside of $\frac{23\times22}{2}$ substitutions which are of the form 
$$(z,az+b),$$
being a quadratic residue of $23$; and if we notice that the substitutions of that function who have the fewest letters are those which have $16$, we see that all the derivatives of a circular substitution of $23$ letters and of a substitution of less than $16$ letters taken among these last form the $23\!$ possible substitutions, or only half, depending on whether the second substitution is even or odd.

\section*{Form of the substitution which characterizes all the preceding transitive functions.}
Let us find out in what form we can write the substitution ($D$) of $n^\circ 1$, which is 
$$\text{(D)} \left(x_z',x_{mz}')(x_z'',x_{mz}''\right),$$ 
by asking
$$\gamma^u\equiv{m}, \quad(\bmod .\frac{p-1}{2}),$$
when we restore the $x$ without accents: $x_1,x_2,\dots$.

The first part\\
$$\text{(L)} \left(x_z',x_{mz}'\right)$$ 
of the substitution (D) is carried out on the quantities $x_0',x_1',x_2',\dots$, or
$$x_{g^0},x_{g^2},x_{g^4},\dots,x_{g^{p-3}},$$
whose indices are quadratic residuals according to the modulus $p$; it comes down to this substitution on the exponents of the even powers of $g$
$$
(2 t, 2 m t), \quad(\bmod . p-1);
$$
thus $g^{2 t}$ is changed to $g^{2 m t}$; therefore by denoting $a$ an arbitrary quadratic residue with respect to $p$, the substitution $(\mathrm{L})$ can be written.
$$
\left(a, a^{m}\right).
$$

Consider the second part

$$
\text{(M)}\qquad \left(x_{z}^{\prime \prime}, x_{m_{z}}^{\prime \prime}\right) 
$$
substitution (D), it is performed on $x_{g}, x_{g^{3}}, x_{g^5}, \dots ;$ only one of these quantities is not moved, and we represent it by $x_{g^{\beta}}$, $\beta$ being odd. So in general, $g^{\beta+2 t}$ will be changed to $g^{\beta+2 m t}$. If we represent the non-residues of $p$ through $b$ and we ask, depending on the module $p$, 
$$g^{\beta+2t}\equiv b,$$
we’ll have
$$
g^{\beta+2 m t}=g^{(1-m) \beta} b^{m}
$$
and the substitution (M) can be written as
$$
\left(b, g^{(1-m) \beta} b^{m}\right)
$$

We have just obtained two formulas to represent the substitution ($D$), depending on whether they are indices which are quadratic residues or indices which are non-residual; but we can represent it by the only formula

$$\left(z, \mathrm{A}z^{\frac{p-1}{2}+m}+\mathrm{B}z^m\right)$$
In fact, if $a$ is a residue of $p$, we have $a^{\frac{p-1}{2}}\equiv 1(\bmod . p)$, and, by expressing that the substitution changes $a$ to $a^{m}$, we get
$$
\mathrm{A}+\mathrm{B} \equiv 1.
$$
If $b$ is a non-residue of $p$, we have $b^{\frac{p-1}{2}}\equiv{-1}(\bmod . p)$, and, in using the second part of the substitution, we have
$$-\mathrm{A}+\mathrm{B} \equiv g^{(1-m) \beta}.$$
We conclude
$$
\mathrm{A} \equiv \frac{1-g^{\beta(1-m)}}{2  },\quad \mathrm{B} \equiv \frac{1+g^{\beta(1-m)}}{2}.
$$
The substitutions that we have designated, in the $n^{\circ}\mathbf{2}$, by $\mathbf{1}^{\circ}, \mathbf{2}^{\circ}, \mathbf{3}^{\circ}$  can be written respectively
$$
\begin{array}{l}
\left(z,-2 z^{5}+3 z^{2}\right) (\bmod 7), \\
\left(z, z^{5}\right), \\
\left(z,-z^{5}-2 z^{2}\right);
\end{array}
$$
therefore the two-time transitive function of seven letters is characterized by the substitution $(z, z+1)$ and any of the three preceding substitutions.

The substitutions that we have designated, in the $n^{\circ} \mathbf{3}$, by $2^{\circ}$ and $4^{\circ}$ can be written, according to modulo $11$,
$$
\begin{array}{l}
\left(z, 5 z^{9}-4 z^{4}\right), \\
\left(z, 3 z^{9}-2 z^{4}\right);
\end{array}
$$
therefore the eleven letter twice transitive function is characterized by $(z, z + 1)$ and one of these two substitutions.

We then see that the four times transitive function of $11$ letters is determined by $(z, z+1)$ and one of the two substitutions
$$
\begin{array}{c}
\left(z,-3 z^{7}+4 z^{2}\right) \\
\left(z, 2 z^{7}-z^{2}\right)
\end{array}
$$
Finally the four times transitive function of $23$ quantities is given by $(z, z+1)$ and by substitution
$$
\left(z,-3 z^{15}+4 z^{4}\right) \text {. }
$$

\section*{Reflections on the transitive functions of a prime number of quantities.}


According to the previous calculations, there is a transitive function of $7$ quantities which is twice, two transitive functions of $11$ quantities, one of which is twice and the other four times, finally a single transitive function of $23$ quantities which is four times. If we generally denote by $p$ the number of quantities enclosed in these functions, these functions are invariant by the substitutions of the form

$$
\text{(g)} \left(z, a^{2} z+b\right),\hfill 
$$
and, $\infty $ being the index of a new quantity, they give functions of $p+1$  quantities which are transitive once more and invariant by the substitutions

$$\left(z, \frac{\mathrm{A}z+\mathrm{B}}{\mathrm{C}z+\mathrm{D}}\right)$$
for which $AD-BC$ is a quadratic residue.

Between $11$ and $23$ are the three prime numbers $13,17,$  and $19$, and I had to wonder if there are any transitive functions of $13,17$ or $19$ letters invariant by all the substitutions ($g$). After doing the necessary calculations, I recognized that there were none. From the tests that I have made, I am, moreover, extremely inclined to believe:
that there are no multiple transitive functions of $13$ or of $19$ letters. As for the transitive functions of $17$ quantities, there are three, and which are three times transitive, provided by this theorem that I gave in my Memoir on the functions of several quantities (Chapter II): "If  $\tau$  is  a  divisor  of $\nu$, there exists a three times transitive function of $p^{\nu}+1$ quantities, which has $\frac{1.2 .3 \ldots\left(p^{\nu}-2\right) \times \tau}{\nu}$  values. We get them by making $p=2, \nu=4$ and $\tau$ successively equal to $1,2,4.$"

Thus the several times transitive functions of $7, 11$ and $23$ quantities, and those of $8, 12$ and $24$, which have been calculated in the current Memoir are indeed due to the fact that the prime numbers $7, 11$ and $23$ are doubles of numbers first increases by one unit. Moreover, as it is easy to see, a transitive function whose prime number increases by one unit is at least twice transitive.



\section*{Acknowledgements} We are very grateful for the encouragement and support of the parents of our talented young translators many of whom are not yet even in middle school. This was done as part of our project to express $M_{23}$ as additive functions on the finite field of $2^{11}$ elements. 

We would also like to thank the help from two really useful websites. We could never have read some of the characters without the help of this OCR. Also to aid in translating, we have used Google Translate. \\\\
\underline{https://www.i2ocr.com/free-online-french-ocr}\\
\underline{https://translate.google.com/?sl=fr\&tl=en\&op=translate}






\printbibliography




\end{document}